\documentclass[12pt]{article}

\usepackage{makeidx}         
\usepackage{graphicx} 
\makeindex   
\usepackage{amsmath,amsthm}
\usepackage{amssymb,latexsym}
\usepackage{bussproofs}
\usepackage[all]{xy}
\xyoption{2cell}
\xyoption{curve}
\UseTwocells
\input{diagxy}
\CompileMatrices      
\usepackage{tikz}
\usepackage{pdfpages}

\newcommand{\T}{\ensuremath{\mathbb{T}}}
\newcommand{\CC}{\ensuremath{\mathcal{C}}}
\newcommand{\EE}{\ensuremath{\mathcal{E}}}
\newcommand{\BB}{\ensuremath{\mathcal{B}}}

\newcommand{\Set}{\ensuremath{\mathsf{Set}}}

\renewcommand{\hom}{\ensuremath{\mathrm{Hom}}}

\newcommand{\mono}{\ensuremath{\rightarrowtail}}

\newtheorem{theorem}{Theorem}
\newtheorem*{theorem*}{Theorem}
\newtheorem{proposition}[theorem]{Proposition} 
\newtheorem*{proposition*}{Proposition} 
\newtheorem{lemma}[theorem]{Lemma}
\newtheorem*{lemma*}{Lemma}
\newtheorem{corollary}[theorem]{Corollary} 
\newtheorem*{corollary*}{Corollary} 

\theoremstyle{remark}
 
\newtheorem*{remarks*}{Remarks}

\theoremstyle{definition}
\newtheorem{definition}[theorem]{Definition}
\newtheorem*{definition*}{Definition}

\newcommand{\myemph}[1]{\emph{#1}}

\begin{document}

\title{
Sheaf Representations and Duality in~Logic
}
\author{
Steve Awodey\thanks{
The author is grateful for the support of the Centre for Advanced Study (CAS) in Oslo, Norway, which funded and hosted the research project ``Homotopy Type Theory and Univalent Foundations'' during the 2018/19 academic year, as well as that of the Air Force Office of Scientific Research through MURI grant FA9550-15-1-0053.}\\
Carnegie Mellon University
}

\date{\today}

\maketitle

\centerline{\emph{In memory of Joachim Lambek}}
\bigskip

\begin{abstract}
The fundamental duality theories relating algebra and geometry that were discovered in the mid-20th century can also be applied to logic via its algebraization under categorical logic.  They thereby result in known and new completeness theorems.  This idea can be taken even further via what is sometimes called ``categorification'' to establish a new connection between logic and geometry, a glimpse of which can also be had in topos theory.  
\end{abstract}

\section*{Preface}

Shortly after finishing my PhD thesis, I received a friendly letter from Professor Lambek in which he expressed interest in a result of mine that extended his work with Moerdijk \cite{LM}. He later cited my result in some papers on the philosophy of mathematics (including \cite{L1,L2}), in which he developed a congenial position that attempted to reconcile the various competing ones in foundations on the basis of results concerning the free topos, the sheaf representations considered here, and related considerations from categorical logic.    

The particular result in question, discussed in section \ref{section:localtopos} below, extends prior results by Lambek and Moerdijk \cite{LM} and Lambek \cite{L2}, and was later extended further in joint work with my  PhD students, first Henrik Forssell \cite{AF, For}, and then Spencer Breiner \cite{AB,B}.  This line of thought is, however, connected to a deeper one in modern mathematics, as I originally learned from the papers of Lambek.  That insight inspired  my  original contribution and also the later joint work with my students, and it continues to fascinate and inspire me.   The purpose of this survey is to sketch that line of thought, which owes more to Lambek than to anyone else.  

The main idea, in a nutshell, is that the ground-breaking duality theories developed in the mid-20th century can also be applied to logic, via its algebraization under categorical logic, and they thereby result in known and new completeness theorems.  This insight, which is already quite remarkable, can as it turns out be taken even further---via what is sometimes called ``categorification''---to establish an even deeper relation between logic and geometry, a glimpse of which can also be had in topos theory, and elsewhere.

\section{Gelfand duality}

Perhaps the ur-example of the sort of duality theory that we have in mind is the relation between topological spaces and commutative rings given by Gelfand duality (see \cite{J} Ch.~4).  To give a brief (and ahistorical) sketch, let $X$ be a  space and consider the ring of real-valued continuous functions on $X$, with pointwise algebraic operations,
\[
 \CC(X)\ =\ \mathsf{Top}(X, \mathbb{R}).
 \]
 This construction is a (contravariant) functor from ``geometry'' to ``algebra'',
 \[
 \mathcal{C} : \mathsf{Top}^\mathsf{op} \to \mathsf{CRng}.
 \]
The functor $ \mathcal{C}$ can be shown  to be full and faithful if we restrict to compact Hausdorff spaces $X$ and (necessarily bounded) continuous functions  $\mathcal{C}^*(X)$,
 \[
 \mathcal{C}^* : \mathsf{KHaus}^\mathsf{op} \hookrightarrow \mathsf{CRng}.
 \]
It then requires some further work to determine exactly \emph{which} commutative rings are of the form $\mathcal{C}^*(X)$ for some space $X$.   These are called (commutative) \emph{$C^*$-algebras}, and they can be characterized as commutative rings $A$ satisfying the following conditions (\cite{J}, \S4.4):
\begin{enumerate}
\item the additive group of $A$ is divisible and torsion free,
\item $A$ has a partial order compatible with the ring structure and such that $a^2 \geq 0$ for all $a\in A$,
\item $A$ is Archemedian, i.e.\ for every  $a\in A$ there is an integer $n$ such that $n\cdot 1_A \geq a$,
\item if $1_A\geq n\cdot a$ for all positive integers $n$, then $a\leq 0$,
\item $A$ is complete in the norm $$||a|| = \inf\{ q\in \mathbb{Q}^+\ |\ q\cdot 1_A \geq a\text{ and } q\cdot 1_A \geq -a\}\,.$$
\end{enumerate}
There are many equivalent specifications (most using complex numbers in place of reals).   Henceforth all rings are assumed to be commutative with unit.

 \begin{theorem}[Gelfand duality]
 The category $\mathsf{KHaus}$ of compact Hausdorff spaces is dual to the category $\mathsf{C^*Alg}$ of $C^*$-algebras and their homomorphisms, via the functor $\mathcal{C}^*$:
 \[
 \mathsf{KHaus}^\mathsf{op} \simeq \mathsf{C^*Alg}.
 \]
 \end{theorem}

How can we recover the space $X$ from its ring of functions $\mathcal{C}^*(X)$?

\begin{itemize}
\item The points  $x\in X$ determine maximal ideals in the ring $\mathcal{C}^*(X)$,
\[
M_x = \{\, f : X\to \mathbb{R}\ |\ f(x) = 0\, \}\,,
\]
and every maximal ideal in $\mathcal{C}^*(X)$ is of this form for a unique $x\in X$.

\item For any ring $A$, the (Zariski) topology on the set $\mathsf{Max}(A)$ of maximal ideals has a basis of open sets of the form:
\[
B_a\, =\, \{M \in X\ |\ a\notin M\, \}\, ,\qquad a\in A.
\]

\item If $A$ is a \emph{$C^*$-algebra}, then this specification  determines a compact Hausdorf space $X = \mathsf{Max}(A)$ such that $A\cong \mathcal{C}^*(X)$.
\end{itemize}

A key step in the proof is the following:
 \begin{proposition}
 Let $A$ be a $C^*$-algebra.  For any maximal ideal $M$ in $A$, the quotient field 
 $A/M$ is isomorphic to  $\mathbb{R}$. 
    \end{proposition}
It follows that there is an injection of rings,
$$A \mono \prod_{M\in \mathsf{Max}(A)}\!\!A/M\ \cong\ \mathbb{R}^{\mathsf{Max}(A)}\,.$$
The image of this map can be shown to consist of the Zariski continuous functions, i.e.\  it is $\mathcal{C}^*(\mathsf{Max}(A))$.


\section{Grothendieck duality for commutative rings}

Grothendieck extended the Gelfand duality from $C^*$-algebras to \emph{all} commutative rings by generalizing on the ``geometric'' side from compact Hausdorff spaces to the new notion of \emph{(affine) schemes},
 \[
\mathsf{Scheme}_\mathsf{aff}^\mathsf{op}\ \simeq\ \mathsf{CRng}.
 \]
The essential difference is to generalize the ``ring of values'' from the constant ring $\mathbb{R}$ to a ring $\mathcal{R}$ that ``varies continuously over the space $X$'', i.e.\ a \emph{sheaf of rings}.   
The various rings $\mathcal{R}_x$ that are the stalks of $\mathcal{R}$ generalize the \emph{local rings} of real-valued functions that vanish at the points $x\in X$.
This change allows \emph{every} commutative ring $A$ to be seen as a ring of continuous functions on a suitable space $X_A$ (the prime spectrum of $A$), where the values of the functions are in a suitable sheaf of (local) rings $\mathcal{R}$ on $X_A$ (see \cite{J} Ch.~5).


\begin{definition} A  ring is called \emph{local} if it has a unique maximal ideal. \\
Equivalently, if $0\neq 1$, and
\begin{equation}\label{eq:localring}
\text{$x+y$ is a unit}\quad\text{implies}\quad\text{$x$ is a unit or $y$ is a unit}.
\end{equation}
\end{definition}
\begin{theorem}[Grothendieck sheaf representation]
Let $A$ be a ring.  There is a space $X_A$ with a sheaf of rings $\mathcal{R}$ such that:
\begin{enumerate}
\item for every $p\in X_A$, the stalk $\mathcal{R}_p$ is a local ring, 
\item there is an isomorphism, $$A\cong \Gamma(\mathcal{R})\,,$$
where $\Gamma(\mathcal{R})$ is the ring of global sections.
\end{enumerate}
Thus every ring is isomorphic to the ring of global sections of a sheaf of local rings.
\end{theorem}


The \myemph{space} $X_A$ in the theorem is the \emph{prime spectrum} $\mathsf{Spec}(A)$ of the ring~$A$:
 \begin{itemize}
\item points $p\in \mathsf{Spec}(A)$ are prime ideals $p\subseteq A$,
\item the (Zariski) topology has basic opens of the form:
 $$B_f = \{ p\in \mathsf{Spec}(A)\ |\ f\not\in p \}, \quad f\in A.$$
\end{itemize}
Note the similarity to the space $\mathsf{Max}(A)$ of maximal ideals from the Gelfand case.  
Unlike that case, however, the functor 
\[
\mathsf{Spec} : \mathsf{CRng}^\mathsf{op} \to \mathsf{Top}
 \]
is not full, and so we need to equip the spaces  $\mathsf{Spec}(A)$ with an additional structure.

The \myemph{structure sheaf} $\mathcal {R}$ is determined at a basic open set $B_f$ by ``localizing'' $A$ at $f$,
\[
\mathcal{R}(B_f) = [f]^{-1}A
\]
where $A \rightarrow [f]^{-1}A$ freely inverts all of the elements $f, f^2, f^3, \dots$.
\medskip

The \myemph{stalk} $\mathcal {R}_p$ of this sheaf at a point $p\in\mathsf{Spec}(A)$ is then seen to be the localization of $A$ at $S_p = A\setminus p$,
\[
\mathcal{R}_p = S^{-1}_p{A}\,.
\]


The \myemph{affine scheme} $(\mathsf{Spec}(A), \mathcal {R})$ presents $A$ as a ``ring of continuous functions'' in the following sense: 
\begin{itemize}
\item each element $f\in A$ determines a ``continuous function'',
\[
\hat{f} : \mathsf{Spec}(A) \to \mathcal {R}\,,
\]
except that the ring $\mathcal{R}$ is itself ``varying continuously over the space $\mathsf{Spec}(A)$'' -- i.e.\ it is a sheaf -- and the function $\hat{f}$ is then a global section of the sheaf $\mathcal{R}$.  

\item Each stalk $\mathcal{R}_p$ is a local ring, with a unique maximal ideal, corresponding to ``those functions  $\hat{f} : \mathsf{Spec}(A) \to \mathcal{R}$ that vanish at $p$''.

\item $(\mathsf{Spec}(A), \mathcal {R})$ is a ``representation'' of $A$ in the sense that $f\mapsto\hat{f}$ is an isomorphism of rings
\[
A \cong \Gamma(\mathcal{R})\,.
\]
\end{itemize}

There is always an injective homomorphism from the global sections of a sheaf into the product of all the stalks,
\[
\Gamma(\mathcal{R}) \rightarrowtail \prod_{p}\mathcal{R}_p \,.
\]
Thus we have the following:

\begin{corollary}[``Subdirect-product representation'']
Every ring $A$ is isomorphic to a {sub}ring of a ``direct product'' of local rings.
I.e.\ there is an injective ring homomorphism
\[
A \rightarrowtail \prod_{p}\mathcal{R}_p \,,
\]
where the $\mathcal{R}_p$ are all local rings.
\end{corollary}

\section{Lambek-Moerdijk sheaf representation for toposes}

\begin{definition} Call a  (small, elementary)  topos $\EE$  \emph{sublocal}\footnote{
In \cite{LM}, and elsewhere, the term \emph{local} was used for the concept here called \emph{sublocal}, and another term was then required in \cite{A} for the stronger condition that we now call \emph{local} in Definition \ref{def:localtopos} below.
} 
if its subterminal lattice $\mathsf{Sub}_\EE(1)$ has a unique maximal ideal.
Equivalently, $0\ncong 1$ and for $x,y\in \mathsf{Sub}_\EE(1)$: 
\[
x\vee y = 1\quad\text{implies}\quad x=1\ \text{or}\ y=1\,.
\]
\end{definition}
Note the formal analogy to the concept of local ring.
In \cite{LM} the following analogue of the Grothendieck sheaf representatation for rings is given for toposes (henceforth, \emph{topos} unqualified will mean small, elementary topos):
\begin{theorem}[Lambek-Moerdijk sheaf representation]
Let $\EE$ be a topos.  There is a space $X_\EE$ with a sheaf of toposes $\widetilde\EE$ such that:
\begin{enumerate}
\item for every $p\in X_\EE$, the stalk $\widetilde\EE_p$ is a sublocal topos, 
\item for the topos $\Gamma(\widetilde\EE)$ of global sections, there is an isomorphism, $$\EE\cong \Gamma(\widetilde\EE)\,.$$
\end{enumerate}
Thus every topos is isomorphic to the topos of global sections of a sheaf of sublocal toposes.
\end{theorem}


The space $X$ mentioned in the theorem is what may be called the \emph{subspectrum of the topos},  $X=\mathsf{sSpec}(\EE)$;   
it is the prime ideal spectrum of the distributive lattice $\mathsf{Sub}(1)$:
\begin{itemize}
\item the points $P\in \mathsf{sSpec}(\EE)$ are prime ideals $P\subseteq \mathsf{Sub}(1)$,
\item the topology has basic opens of the form:
$$B_q = \{ P\in \mathsf{Spec}(\EE)\ |\ q\not\in P \}, \quad q\in\mathsf{Sub}(1)\,.$$
\end{itemize}
Note the close analogy to the space $\mathsf{Spec}(A)$ for a commutative ring $A$.

The lattice of all open sets of $\mathsf{sSpec}(\EE)$ is then (isomorphic to) the ideal completion of the lattice $\mathsf{Sub}(1)$,
$$\mathcal{O}(\mathsf{Spec}(\EE)) = \mathsf{Idl}(\mathsf{Sub}(1))\,.$$

%


Next, let us define a \emph{structure sheaf} $\widetilde{\EE}$ on $\mathsf{sSpec}(\EE)$ by ``slicing'' $\EE$ at $q \in\mathsf{Sub}(1)$,
\[
\widetilde{\EE}(B_q) = \EE/q\,.
\]
This takes the place of the localization of a ring $A$ at a basic open $B_f$:
\[
\mathcal{R}_A(B_f) = [f]^{-1}A\,.
\]
Note that $\EE/q$ ``inverts'' all those elements $p\in \mathsf{Sub}(1)$ with $q\leq p$, in the sense that the canonical map 
$q^* : \EE \to \EE/q$ takes every $p \mono 1$ to $q\wedge p \mono q$, and so if $q\leq p$ then $q^*p = 1_q : q \mono q$.

The fact that $\widetilde{\EE}$ is indeed a sheaf on $\mathsf{sSpec}(\EE)$ comes down to showing that, for any $p, q \in \mathsf{Sub}(1)$, there is a canonical equalizer of toposes (and logical morphisms),
\[
\EE/p\vee q\ \mono\ \EE/p \times \EE/q\ \rightrightarrows\ \EE/p\wedge q \,.
\]
This in turn says that in a diagram of the form: 
 
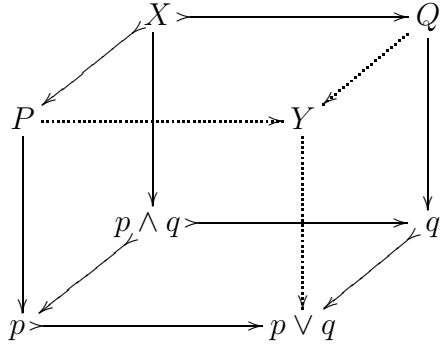
\[ 
 \xymatrix{ 
& \ \ X\ \ar@{>->}[ld] \ar[dd] \ar@{>->}[rr] && Q \ar@{.>}[ld] \ar[dd] & \\
 P \ar[dd] \ar@{.>}[rr] && Y \ar@{.>}[dd]\\ 
 & p \wedge q\ \ar@{>->}[ld] \ar@{>->}[rr] && \ q \ar@{>->}[ld] & \\
  p\ \ar@{>->}[rr] && p\vee q} 
 \] 
 with a pushout-pullback of monos in the base, and the two vertical squares involving $X$ given as pullbacks, one can complete the cube as indicating by first forming the pushout $Y$ on the top face, and then obtaining the front vertical map from $Y$, and the resulting new vertical faces will then also be pullbacks. This is a rather special ``descent condition'' for the presheaf $\EE/-$.
 
The stalk $\widetilde{\EE}_P$ of this sheaf at a point $P\in \mathsf{sSpec}(\EE)$ is computed as the filter-quotient of $\EE$ over the complement of the prime ideal $P\subseteq  \mathsf{Sub}_{\EE}(1)$, i.e.\ the prime filter $P^c = \mathsf{Sub}(1)\!\setminus\! P$.  Thus for the stalk we have the (filtered) colimit (taken in $\mathsf{Cat}$, but again a topos):
\[
\widetilde{\EE}_P = \varinjlim_{q\not\in P} \EE/q\,.
\]
For this stalk topos, one then has the subterminal lattice:
\[
\mathsf{Sub}_{\widetilde{\EE}_P}(1)\cong \mathsf{Sub}_{\EE}(1)/P^c\,,
\]
where $\mathsf{Sub}_{\EE}(1)/P^c$ is the quotient Heyting algebra by the prime filter $P^c$.  Since for the prime filter $P^c$ we have $p \vee q \in P^c$ implies $p \in P^c$ or $q \in P^c$, it thus follows that the stalk topos $\widetilde{\EE}_P$ is indeed sublocal.

Finally, for the global sections of $\widetilde{\EE}$ we  have simply:
$$\Gamma(\widetilde{\EE}) \cong \widetilde{\EE}(B_\top) = \EE/1 \cong \EE\,.$$
Thus the topos of global sections of $\widetilde{\EE}$ is indeed isomorphic to the original topos $\EE$.  In this way, $\EE$ is isomorphic to the global sections of a sheaf of sublocal toposes.
\medskip

Again, there is always an injection from the global sections into the product of the stalks, which in this case gives a conservative logical morphism of the form
\[
\EE \cong \Gamma(\widetilde{\EE}) \rightarrowtail \prod_{P\in \mathsf{sSpec}(\EE)}\!\widetilde{\EE}_P\,.
\]

\begin{corollary}\label{cor:SDPtopos}
Every topos has a conservative logical morphism into a product of sublocal toposes.
\end{corollary}

\subsection{Lambek's modified sheaf representation for toposes}

Now consider  the following  logical interpretation of the sheaf representation theorem for toposes and its corollary.
\begin{itemize}
\item A topos $\EE$ can be regarded as the syntactic category $\EE_\mathbb{T}$ of a theory $\mathbb{T}$ in Intuitionistic Higher-Order Logic (IHOL).  Thus for any sentence $\phi$ in the language of the theory $\mathbb{T}$,
\[ 
\EE_\mathbb{T}\models \phi \qquad\text{iff}\qquad \mathbb{T} \vdash \phi\,.
\]

\item A logical functor $\EE_\mathbb{T}\rightarrow\mathcal{F}$ between toposes induces an interpretation of the language of $\mathbb{T}$ into $\mathcal{F}$; one which, moreover, satisfies all the sentences that hold in $\EE_\mathbb{T}$, and thus a $\mathbb{T}$-model.  Such a functor is called \emph{conservative} if it reflects isomorphisms; for logical functors, this is the same as being faithful.  A conservative logical functor $f : \EE_\mathbb{T}\rightarrowtail\mathcal{F}$ therefore reflects satisfaction, in the sense that for any sentence $\phi$ in the language of $\mathbb{T}$,
\[ 
\mathcal{F} \models \phi^f \qquad\text{implies}\qquad \EE_\mathbb{T}\models \phi\,.
\]
where $\phi^f $ is the interpretation of $\phi$ under the model induced by $f$.

\item A sublocal topos $\mathcal{S}$ is one that is \emph{consistent} $\mathcal{S}\nvDash\bot$ and has the \emph{disjunction property}  
\[
\mathcal{S}\models \phi\vee \psi \qquad\text{iff}\qquad \mathcal{S}\models \phi\  \ \text{or}\  \ \mathcal{S}\models \psi\,,
\]
for all sentences $\phi, \psi$.  Such sublocal toposes are more $\Set$-like than a general one, and can thus be regarded as suitable semantics for logical theories. 

\item  The ``subdirect-product representation'' given by Corollary \ref{cor:SDPtopos} is a logical completeness theorem with respect to interpretations $\EE_\mathbb{T}\rightarrow\mathcal{S}$ of $\mathbb{T}$ into sublocal toposes~$\mathcal{S}$.  It says that, for any theory $\mathbb{T}$ in IHOL, a sentence $\phi$ is provable, $\mathbb{T}\vdash \phi$,  iff it holds in every interpretation of $\mathbb{T}$ in a sublocal topos~$\mathcal{S}$.  Thus IHOL is complete with respect to models in sublocal toposes.

\item The sheaf representation is a Kripke-style completeness theorem for IHOL, with $\widetilde\EE$ as a ``sheaf of possible worlds'' (see \cite{L1}).
\end{itemize}

Under this interpretation, however, the present sheaf representation is not entirely satisfactory, because we would really like the ``semantic toposes'' $\mathcal{S}$ to be even more $\Set$-like, in addition to being sublocal, by also having the \emph{existence property}:
\[
\mathcal{S}\models (\exists x:A)\varphi(x) \qquad\text{iff}\qquad  \mathcal{S}\models \varphi(a)\text{, for some closed $a:A$}\,.
\] 

\begin{definition}\label{def:localtopos}
A topos $\mathcal{S}$ will be called \emph{local} if the terminal object $1$ is both indecomposable and projective, i.e.\ the global sections functor 
\[
\Gamma = \hom_\mathcal{S}(1, - ) : \mathcal{S} \to \Set
\]
preserves all finite coproducts and epimorphisms. Note that a local topos is exactly one that is consistent and has both the disjunction and existence properties.
\end{definition}

In the paper \cite{L1}, Lambek gave the following improvement over the sublocal sheaf representation:

\begin{theorem}[Lambek sheaf representation]
Let $\EE$ be a topos.  
There is a faithful logical functor $\EE\rightarrowtail\mathcal{F}$ 
and a space $X$ with a sheaf of toposes 
$\widetilde{\mathcal{F}}$ such that:
\begin{enumerate}
\item for every $p\in X$, the stalk $\widetilde{\mathcal{F}}_p$ is a  local topos, 
\item for the global sections of $\widetilde{\mathcal{F}}$ there is an isomorphism $\mathcal{F} \cong \Gamma(\widetilde{\mathcal{F}}) $.
\end{enumerate}
Thus every topos is a  subtopos of one that is isomorphic to the global sections of a sheaf of  local toposes.  
\end{theorem}
\medskip
The proof was inspired by the Henkin completeness theorem for higher-order logic \cite{H}, and first performs a sort of ``Henkinization'' of $\EE$ to get a bigger topos $\EE\mono \mathcal{F}$ with witnesses for all existential quantifiers, in a suitable sense. This result then suffices for a subdirect-product embedding of any topos $\EE$ into a product of \emph{local} toposes, and therefore gives the desired logical completeness of IHOL with respect to such toposes, which are much more $\Set$-like.

\section{Local sheaf representation for toposes}\label{section:localtopos}

The result from \cite{A} mentioned above was this:

\begin{theorem}[Local topos sheaf representation]\label{theorem:local}
Let $\EE$ be a topos.  
There is a space $X_\EE$ with a sheaf of toposes $\widetilde{\EE}$ such that:
\begin{enumerate}
\item for every $p\in X_\EE$, the stalk $\widetilde{\EE}_p$ is a local topos, 
\item for the global sections of $\widetilde{\EE}$  there is an equivalence $\EE\simeq \Gamma(\widetilde{\EE})$.
\end{enumerate}
Thus every topos is equivalent to the global sections of a sheaf of local toposes.  
\end{theorem}

As before, this gives a \emph{subdirect-product representation} of $\EE$,
\[
\EE \rightarrowtail \prod_{p\,\in X}\mathcal{S}_p\,,
\]
into a product of local toposes $\mathcal{S}_p = \widetilde{\EE}_p$, and therefore implies the desired logical completeness of IHOL with respect to local toposes.  
This stronger result also gives better ``Kripke semantics'' for IHOL, since the ``sheaf of possible worlds'' (in the sense of \cite{L1}) now has local stalks.
\medskip

For classical higher-order logic, something more can be said:

\begin{lemma}
Every Boolean, local topos $\mathcal{S}$ is well-pointed, i.e.\ the global sections functor,
\[
\Gamma = \hom_\mathcal{S}(1, - ) : \mathcal{S} \to \Set
\]
is faithful.
\end{lemma}

\begin{corollary}
Every Boolean topos is isomorphic to the global sections of a sheaf of well-pointed toposes.  
\end{corollary}
For Boolean toposes $\mathcal{B}$, we therefore have an embedding, 
\[
\mathcal{B} \rightarrowtail \prod_{p\,\in X}\mathcal{S}_p 
\]
as a subdirect-product of \emph{well-pointed} toposes $\mathcal{S}_p$ (this is \cite{F}, Thm 3.22).  The logical counterpart now says:
\begin{corollary}
Classical HOL is complete with respect to models in well-pointed toposes.
\end{corollary}
A well-pointed topos is essentially a model of classical Zermelo set theory (\cite{MM}, {\S}VI.10).  Indeed, it is worth emphasizing that the models of HOL here are \emph{standard} models of classical HOL (i.e.\ with full function and power sets), taken in \emph{varying} models $\mathcal{S}$ of set theory.

Finally, taking the global sections $\Gamma:\mathcal{S}_p \rightarrowtail \Set$ of each well-pointed topos $\mathcal{S}_p$, we get a faithful functor from any Boolean topos $\mathcal{B}$ into a power of~$\Set$:
\[
\mathcal{B} \rightarrowtail \prod_{p\,\in X}\mathcal{S}_p \rightarrowtail \prod_{p\,\in X}\Set \cong \Set^X\,.
\]
However, the various composites $\mathcal{B} \rightarrow \mathcal{S}_p \rightarrowtail \Set$ are now not logical functors, because they need not preserve exponentials; they do, however, preserve the first-order logical structure (they are also exact; thus we have another proof of \cite{F} theorem 3.24).  These composites are exactly what the logician calls a ``Henkin'' or ``non-standard'' model of HOL in $\Set$.  In this way, we recover the familiar ``Henkin completeness theorem for HOL'' \cite{H}:

\begin{corollary}
Classical HOL is complete with respect to Henkin models in~$\Set$.
\end{corollary}

For the proof of the local topos sheaf representation theorem, these ``Henkin models'' will be taken as the points of the space $X_\EE$, which we call the \emph{space of models} (following \cite{BM}).
In the  sublocal case, the points were the \emph{prime ideals} $P\subseteq\mathsf{Sub}(1)$.
These correspond exactly to the \emph{lattice homomorphisms} $$p: \mathsf{Sub}_{\EE}(1)\to \mathbf{2}\,.$$

For the local case, we instead take \emph{coherent functors} $$P: \EE\to\Set\,,$$
which correspond to (Henkin) models of the ``theory'' $\EE$.\footnote{
Of course, the collection of all such functors may be too big to form a \emph{set}.  The remedy, as explained in the paper \cite{A}, is to choose a suitable cardinal bound $\kappa$ on the size of the models $P$.}

The topology on $X_\EE$ can be described roughly as follows (see \cite{A} for more details, but the idea for this topology originates with \cite{JM,BM}; it was also used in \cite{AF,B}).  To simplify things, let us regard \EE\ as a classifying topos for a theory $\T$, and say that a model $P: \EE\to \Set$ ``satisfies'' a sentence $\phi$, which we may identify with its interpretation $\phi \mono 1_\EE$, if $[\![\phi]\!]^P = P(\phi) = 1$. Then we could mimic the subspectrum by taking as a basic open set all those models $P$ that satisfy some fixed $\phi$:  
\[
V_\phi = \{ P\in X_\EE\ |\ P\models \phi \}\,.
\]
However, it turns out that there are too few such basic opens; thus we will also use formulas $\phi(x)$ with free variables .  In order to say when $P\models \phi(x)$ we therefore equip each model $P$ with a ``labelling'' $\alpha : \kappa \rightarrow |P|$ by elements of some fixed, large set $\kappa$, and we then define the notion of satisfaction of a formula by such a labelled model $(P, \alpha) \models \phi(x)$, which we write suggestively as $P \models \phi(\alpha)$.  Thus the points of $X_\EE$ are actually pairs $(P, \alpha)$, and the basic open sets then have the form 
\[
V_{\phi(x)} = \{ (P,\alpha)\in X_\EE\ |\ P\models \phi(\alpha) \}
\]
for all formulas $\phi(x)$. (This description is not entirely accurate, but it gives the idea for present purposes; see \cite{A,AF,B} for details.)

The \emph{structure sheaf} $\widetilde\EE$ on $X_\EE$  is again defined by ``slicing'' $\EE$, 
\[
\widetilde\EE(A)\ =\ \EE/A\,\qquad \text{for $A\in\EE$}\,,
\] 
but now it is first shown to be a \emph{stack} on $\EE$ itself (with respect to the coherent topology). What this means is:
\begin{enumerate}
\item for any $A,B\in \EE$,  the canonical map is an equivalence,
\[
\EE/A+B\ \simeq\ \EE/A \times \EE/B\,,\qquad 
\]

\item for any epimorphism $e:B\epi A$, the canonical map is an equivalence,
\[
\EE/A\ \simeq\ \mathrm{des}(\EE/B, e)\,,
\]
where $\mathrm{des}(\EE/B, e)$ is the category of objects of $\EE/B$ equipped with \emph{descent data} with respect to $e:B\to A$.
\end{enumerate}

The stack is then strictified to a \emph{sheaf} of categories (see \cite{A}), and then finally transferred from $\EE$ to the space $X_\EE$ of models using a topos-theoretic covering theorem due to Butz and Moerdijk \cite{BM}.  Call the resulting sheaf of categories on $X_\EE$ again $\widetilde\EE$.

The \emph{stalk} $\widetilde{\EE}_{(P,\alpha)}$ of the (transferred) sheaf at a point $(P,\alpha)$ can be calculated as the colimit,
\[
\widetilde{\EE}_{(P,\alpha)}\ =\ \varinjlim_{A\,\in\int\!{P}} \EE/A, 
\]
where the (filtered!) category of elements $\int_\EE\!{P}$ of the model $P:\EE\rightarrow\Set$ takes the place of the prime filter.   As a key step, one shows that these stalks are indeed \emph{local} toposes whenever $P:\EE\rightarrow\Set$ is a coherent functor.

Finally, for the \emph{global sections} functor $\Gamma:\mathsf{Sh}(X_\EE)\rightarrow\Set$, we still have:
$$\Gamma(\widetilde{\EE})\ \simeq\ \EE/1\ \cong\ \EE\,.$$  In this way, $\EE$ is indeed equivalent to the topos of global sections of a sheaf of local toposes on a space.



\section{Stone duality for Boolean algebras}

The foregoing sheaf representations for toposes suggest an analogous treatment for \emph{pretoposes}, which would actually be somewhat better, because the $\Set$-valued models used for the points (and coming from the global sections of the stalks) would then all be \emph{standard} models, rather than Henkin style, non-standard models. This suggests the possibility of a \emph{duality theory for first-order logic}, analogous to that for affine schemes and commutative rings, with  pretoposes playing the role of rings, the space of models playing the role of the prime spectrum, and the sheaf representation providing a structure sheaf.

This is more than just an analogy: it is a generalization of the classical Stone duality for Boolean algebras (= Boolean rings). From a logical point of view, the classical duality theory for Boolean algebras is the propositional case of the first-order one that we are proposing for pretoposes.  (There is also a generalization from classical to intuitionistic logic, which is less of a stretch.)  Thus let us briefly review the ``propositional case'' of classical Stone duality for Boolean algebras, before proceeding to the ``first-order'' case of pretoposes.  

Recall (e.g.\ from \cite{J}, Ch.~5) that for a Boolean algebra $B$ we have the Stone space $\mathsf{Stone}(B)$, which is defined exactly as was the subterminal lattice $\mathsf{Sub}_{\EE}(1)$ of a topos $\EE$, i.e. $\mathsf{Stone}(B)=\mathsf{Spec}(B)$ is the prime spectrum of $B$ (prime ideals in a Boolean algebra are always maximal, thus are exactly the complements of the ultrafilters, which are the usual points of $\mathsf{Stone}(B)$).  We can represent the \emph{points} $p\in \mathsf{Spec}(B)$ as Boolean homomorphisms,
\[
p : B\to \mathbf{2}\,.
\]
And we can recover the Boolean algebra $B$ from the space $\mathsf{Spec}(B)$ as the \emph{clopen subsets}, which are represented by continuous maps,
\[
f : \mathsf{Spec}(B)\to \mathbf{2}\,,
\]
where (the underlying set of) $\mathbf{2}$ is given the discrete topology. Note that this is also a sheaf representation -- but a constant one!  The stalks are  \emph{local} Boolean algebras, which are always just~$\mathbf{2}$.  

\emph{Stone's representation theorem} for Boolean algebras then says that there is always an injective homomorphism,
\[
B \mono \mathbf{2}^X \cong \mathcal{P}(X)\,,
\]
for a set $X$, which we can take to be the set of points of $\mathsf{Spec}(B)$, i.e.\ the ultrafilters.  This is therefore the usual subdirect-product embedding resulting from the sheaf representation.

There is, moreover, a contravariant equivalence of categories,
\[
\xymatrix{ 
\mathsf{Bool}  \ar@/{}^{1pc}/[rr]^{\mathsf{Spec}}     & \quad\simeq &  \mathsf{Stone}^{\mathsf{op}}  \ar@/{}^{1pc}/[ll]^{\mathsf{Clop}} \,. \\
} 
\]
Both of the functors $\mathsf{Spec}$ and $\mathsf{Clop}$ are given by homming into $\mathbf{2}$, albeit in two different categories.

Logically, a Boolean algebra $B$  is always the ``Lindenbaum-Tarski algebra'' of a theory $\mathbb{T}$ in \emph{propositional} logic, and a Boolean homomorphism $B\to \mathbf{2}$ is then the same thing as a $\T$-model, i.e.\ a ``truth-valuation''.  Thus the points of $\mathsf{Spec}(B)$ are \emph{models} of the propositional theory $\mathbb{T}$.
We are going to generalize this situation by replacing Boolean algebras with (Boolean) pretoposes, representing \emph{first-order} logical theories, and replacing $\mathbf{2}$-valued models with $\mathsf{Set}$-valued models.

\section{Stone duality for Boolean pretoposes}

M.~Makkai \cite{M} has discovered a Stone duality for Boolean pretoposes with respect to what he terms \emph{ultragroupoids} on the geometric/semantic side. These are groupoids (of models and isomorphisms) equipped with a primitive structure of ultraproducts of models, together with groupoid homomorphisms that preserve ultraproducts.  The result is an equivalence of categories:

\[
\xymatrix{ 
\mathsf{BoolPreTop}  \ar@/{}^{1pc}/[rr]     & \simeq &  \mathsf{UltraGpd}^{\mathsf{op}}  \ar@/{}^{1pc}/[ll] \\
} 
\]
\smallskip

\noindent which, as in the propositional case, is mediated by homming into a special object, now $\Set$ in place of $\mathbf{2}$.  This replacement, and the remarkable duality theory that results, is an instance of what is sometime called ``categorification'', an idea that plays a guiding role throughout categorical logic.  It follows in particular that every Boolean pretopos $\mathcal{B}$ has a pretopos embedding into a power of $\Set$.
\[
\mathcal{B} \mono \Set^X\,,
\]
where $X$ is a set of ``models'', i.e.\ pretopos functors $M : \mathcal{B}\to \Set$.  

We will show below that this last fact---which is essentially G\"odel's completeness theorem for first-order logic---is also a ``subdirect-product representation'' resulting from a sheaf representation of $\mathcal{B}$.  But first we need to make a suitable ``space of models''.

In joint work with H.~Forssell~\cite{AF,For} Makkai's ultragroupoids of models were replaced by \emph{topological} groupoids of models, equipped with a Stone-Zariski type logical topology similar to the one used above for the local sheaf representation for toposes. In overview, our (topological) generalization of Stone duality from Boolean algebras to Boolean pretoposes works like this:
\medskip

\begin{center}
\begin{tabular}{c|c}
Boolean algebra $B$ & Boolean pretopos $\mathcal{B}$ \\
propositional theory & first-order theory \\
\\
\hline\\
 homomorphism & pretopos functor\\
$B\to \mathbf{2}$ &  $\mathcal{B}\to \mathsf{Set}$ \\
truth-valuation & elementary model \\
\\
\hline\\
 topological space & topological groupoid \\
 $\mathsf{Spec}(B)$ & $\mathsf{Spec}(\mathcal{B})$ \\
of all valuations & of all models and isos \\
\\
\hline\\
 continuous function & continuous functor \\
$\mathsf{Spec}(B) \to \mathbf{2}$  & $\mathsf{Spec}(\mathcal{B}) \to \mathsf{Set}$\\
clopen set & coherent sheaf 
\end{tabular}
\end{center}


To give a bit more detail of a few of the steps:
\begin{itemize}
\item The spectrum $\mathsf{Spec}(\mathcal{B})$ of a Boolean pretopos $\mathcal{B}$ is not just a space, but a \emph{topological groupoid}, consisting of a space of (labelled) models $(M,\alpha)$ and a space of isos $i : M \cong N$.  These are topologized by a \emph{logical topology} of the same kind already considered, where the basic opens (of the space of models) are determined by satisfaction of formulas, $$V_{\phi(x)} = \{ (M,\alpha)\in \mathsf{Spec}(\mathcal{B})\ |\ M\models \phi(\alpha)\}\,.$$

\item Morphisms $f : \mathsf{Spec}(\mathcal{B}) \to \mathsf{Spec}(\mathcal{B'})$ are just continuous groupoid homomorphisms.  Every pretopos functor $F : \mathcal{B'} \to \mathcal{B}$ gives rise to such a homomorphism, essentially by precomposition, since 
$$\mathsf{Spec} : \mathsf{BoolPreTop} \to \mathsf{StoneTopGpd}^{\mathsf{op}}$$ 
is representable,
\[
\mathsf{Spec}(\mathcal{B}) \simeq \mathsf{BoolPreTop}(\mathcal{B}, \Set)\,.
\]
Thinking of such a pretopos functor $F :\mathcal{B'} \to \mathcal{B}$ as a ``translation of theories'', the semantic functor $\mathsf{Spec}(F)$ acts on models in the corresponding way.

\item Recovering $\mathcal{B}$ from $\mathsf{Spec}(\mathcal{B})$ amounts to recovering an elementary theory (up to pretopos completion) from its models.  This is done using hard results from topos theory due mainly to Joyal-Tierney and Joyal-Moerdijk \cite{JT,JM,BM}. Specifically, one shows that the category
of \emph{equivariant} sheaves on the topological groupoid $\mathsf{Spec}(\mathcal{B})$ is equivalent to the (Grothendieck) topos of sheaves on $\mathcal{B}$ for the coherent topology,
\[
\mathsf{Sh}_{\mathrm{eq}}(\mathsf{Spec}(\mathcal{B}))\ \simeq\ \mathsf{Sh}(\mathcal{B})\,.
\]
Logically, this gives two different presentations of the (Grothendieck) classifying topos of a first-order theory $\T$, such that $\mathcal{B} = \mathcal{B}_\T$ is the pretopos completion of (the syntactic category of) $\T$, and  $\mathsf{Spec}(\mathcal{B})$ is then the groupoid of $\T$-models.

It follows that $\mathcal{B}$ is equivalent to the subcategory of \emph{coherent objects} of this topos; thus $\mathcal{B}$ is equivalent to the category of coherent, equivariant sheaves on the topological groupoid $\mathsf{Spec}(\mathcal{B})$.  
These can be shown to correspond to certain continuous homomorphisms $\mathsf{Spec}(\mathcal{B}) \rightarrow \Set$, where the latter is the topological groupoid of sets, equipped with a suitable topology.  In this sense, the coherent, equivariant sheaves generalize the clopen sets in a Stone space.

\end{itemize}

Unlike in the case of Boolean algebras, however, and unlike in Makkai's theorem using ultragroupoids, we do not have an equivalence of categories, but only an \emph{adjunction} \cite{AF,For}:
 
\begin{theorem}[Awodey-Forssell]
There is a contravariant adjunction,
\[
\xymatrix{ 
\mathsf{BoolPreTop}  \ar@/{}^{1pc}/[rr]^{\mathsf{Spec}}     &&  \mathsf{StoneTopGpd}^{\mathsf{op}}  \ar@/{}^{1pc}/[ll]^{\mathsf{Coh}} \,,
} 
\]
in which both functors are given by homming into $\mathsf{Set}$.
\end{theorem}

In particular, the ``semantic'' functor,
\[
\mathsf{Spec} : \mathsf{BPreTop} \longrightarrow \mathsf{StoneTopGpd}^{\mathsf{op}}
\]
 is \emph{not full}: there are continuous functors between the groupoids of models that do not come from a ``translation of theories''.
Compare the case of commutative rings $A, B$, where an arbitrary continuous function $$f : \mathsf{Spec}(B) \to \mathsf{Spec}(A)$$ need not come from a ring homomorphism $h : A\to B$.  

We can of course characterize the ``semantic functors'' arising from a pretopos morphism as those that pull coherent sheaves back to coherent sheaves.  Such ``coherent'' maps  $f : \mathsf{Spec}(\mathcal{B}) \to \mathsf{Spec}(\mathcal{B'})$ will then correspond to pretopos maps $F : \mathcal{B'} \to \mathcal{B}$, simply  by $f(M) \cong M\circ F$. 



\section{Sheaf representation for pretoposes}

We now want to cut down the morphisms between the semantic groupoids $\mathsf{Spec}(\mathcal{B})$ to just the coherent ones that come from pretopos functors.  We will do this by endowing $\mathsf{Spec}(\mathcal{B})$ with additional structure that is preserved by all such ``syntactic'' maps.   Specifically, as for rings and affine schemes, we can equip the spectrum $\mathsf{Spec}(\mathcal{B})$ of the pretopos $\mathcal{B}$ with a ``structure sheaf'' $\widetilde{\mathcal{B}}$,  defined just as in the sheaf representation for toposes: 
\begin{itemize}
\item Start with the pseudofunctor $\widetilde{\mathcal{B}} : \mathcal{B}^{\mathrm{op}}\to\mathsf{Cat}$ with,
\[
\widetilde{\mathcal{B}}(X) \cong \mathcal{B}/X\,,\qquad X\in \BB\,.
\]
The prestack $\widetilde{\mathcal{B}}$ is actually a \emph{stack} for the coherent topology, because $\mathcal{B}$ is a pretopos.
 
\item Strictify $\widetilde{\BB}$ to get a sheaf of categories (also called $\widetilde{\mathcal{B}}$) on $\mathcal{B}$.  The ``stalk'' of $\widetilde{\mathcal{B}}$ at a ``point'' $M : \EE\to \Set$ (a pretopos functor) is then 
\[
\widetilde{\EE}_{M}\ \simeq\ \varinjlim_{A\,\in\int\!{M}}\EE/A, 
\]
which is a \emph{local} Boolean pretopos ($1$ is indecomposable and projective).

\item There is an equivalence of Grothendieck toposes,
\[
\mathsf{Sh}(\BB)\ \simeq\ \mathsf{Sh}_{\mathsf{eq}}(\mathsf{Spec}(\mathcal{B})),
\]
between sheaves on the pretopos $\BB$, for the coherent Grothendieck topology, and equivariant sheaves on the topological groupoid $\mathsf{Spec}(\mathcal{B})$ of (labelled) models. 

\item Move $\widetilde{\mathcal{B}}$ across this equivalence in order to get an equivariant sheaf on $\mathsf{Spec}(\mathcal{B})$.  The result (also called $\widetilde{\mathcal{B}}$) is thus a sheaf of local, Boolean pretoposes on $\mathsf{Spec}(\mathcal{B})$.  

\end{itemize}

And from a logical point of view:
\begin{itemize}
\item $\mathcal{B}= \mathcal{B}_\T$ is the Boolean pretopos completion of (the syntactic category of) a theory $\T$ in (classical) FOL, and  
$\mathsf{Spec}(\mathcal{B})$ is then the groupoid of $\T$-models.

\item $\widetilde{\mathcal{B}}$ is a sheaf of ``local theories''.  The stalk $\widetilde{\mathcal{B}}_M$ at a $\T$-model $M$ is a well-pointed pretopos representing the complete theory of $M$, with parameters for all the elements of $M$ added; it is what the logician calls the ``elementary diagram'' of the model $M$.

\item As before, $\widetilde{\mathcal{B}}$ has global sections $\Gamma(\widetilde{\mathcal{B}}) \simeq \mathcal{B}$.  So the original pretopos $\mathcal{B}_\T$ turns out to be the ``theory of all the $\T$-models''.

\item Since each stalk $\widetilde{\mathcal{B}}_M$ is local, and well-pointed,  the global sections functor $\Gamma_M : \widetilde{\mathcal{B}}_M \rightarrowtail \Set$ is a faithful \emph{pretopos} morphism, i.e.\ a model in $\Set$.  In fact, the model $M:\mathcal{B}\to\Set$  is naturally  isomorphic to the composite:
\[
M : \mathcal{B} \simeq \Gamma(\widetilde{\mathcal{B}}) \to \widetilde{\mathcal{B}}_M \stackrel{\Gamma_M}{\to} \Set\,.
\]
\end{itemize}

In sum, we have the following (see \cite{AB,B}):

\begin{theorem}[Awodey-Breiner]
Let $\BB$ be a Boolean pretopos.  
There is a topological groupoid $G$ with an equivariant sheaf of pretoposes $\widetilde{\BB}$ such that:
\begin{enumerate}
\item for every $g\in G$, the stalk $\widetilde{\BB}_g$ is a well-pointed pretopos, 
\item for the global sections of $\widetilde{\BB}$ there is an equivalence $\BB\simeq\Gamma(\widetilde{\BB})$.
\end{enumerate}
Thus every Boolean pretopos is equivalent to the global sections of a sheaf of well-pointed pretoposes.  
\end{theorem}
\medskip

There is again an analogous result for the general (i.e.\ non-Boolean) case, with local pretoposes in place of well-pointed ones in the stalks.  The associated subdirect-product representation is then the following:

\begin{corollary}\label{bptsdp}
For any pretopos $\EE$, there is a pretopos embedding,
\[
\mathcal{E}\rightarrowtail \prod_{g\in X_\EE}\mathcal{E}_{g} 
\]
with each $\EE_g$ a local pretopos and $X_\EE$ the set of points of the topological groupoid $\mathsf{Spec}(\mathcal{E})$.
If moreover $\BB$ is Boolean, then the local pretoposes $\mathcal{B}_{g}$ are all well-pointed, and  $\BB$ therefore embeds (as a pretopos!) into a power of~$\Set$: 
\[
\mathcal{B}\rightarrowtail \prod_{g\in X_\BB}\mathcal{B}_{g} 
\rightarrowtail \prod_{g\in X_\BB}\Set \simeq \Set^{X_\BB}\,.
\]

\end{corollary}
In logical terms, the last statement is essentially the G\"odel completeness theorem for first-order logic, repackaged.  Of course, the proof made use of the equivalent fact that $\BB$ has enough pretopos functors $M:\BB\to\Set$. 

%
%

\section{Logical schemes}

For a Boolean pretopos $\BB$, call the pair $$(\mathsf{Spec}(\mathcal{B}), \widetilde{\mathcal{B}})$$ just constructed an \emph{affine logical scheme}.
A morphism of affine logical schemes
\[
(f, \widetilde{f}) : (\mathsf{Spec}(\mathcal{A}), \widetilde{\mathcal{A}}) \to (\mathsf{Spec}(\mathcal{B}), \widetilde{\mathcal{B}})
\]
consists of a continuous groupoid homomorphism 
\[
f : \mathsf{Spec}(\mathcal{A}) \to \mathsf{Spec}(\mathcal{B}),
\]
together with a pretopos functor over $\mathsf{Spec}(\mathcal{B})$
\[
\widetilde{f} : \widetilde{\mathcal{B}} \to f_*\widetilde{\mathcal{A}}\,.
\]

\begin{theorem}[Awodey-Breiner]
Every pretopos functor $\mathcal{B} \to \mathcal{A}$ induces a morphism of the associated affine logical schemes $\mathsf{Spec}(\mathcal{A}) \to \mathsf{Spec}(\mathcal{B})$.  Moreover, the  functor
\[
\mathsf{Spec} : \mathsf{BoolPreTop} \longrightarrow \mathsf{LogScheme}^{\mathsf{op}}_{\mathsf{aff}}
\]
is full and faithful: every map of schemes comes from an essentially unique map of pretoposes.
\end{theorem}

\begin{corollary}[First-order logical duality]
There is an equivalence,
\[
\mathsf{BoolPreTop} \ \simeq\ \mathsf{LogScheme}^{\mathsf{op}}_{\mathsf{aff}}\,.
\]
\end{corollary}

The category of Boolean pretoposes is thus dual to the category of affine logical schemes.  We can now start to  ``patch together'' affine pieces of the form $(\mathsf{Spec}(\mathcal{B}),\widetilde{\mathcal{B}})$, in order to make a general notion of a ``logical scheme'', consisting of a topological groupoid of structures not tied to any one theory, equipped with a sheaf of local theories, and locally equivalent to an affine scheme.  The first few steps in this direction are explored in \cite{B}.




\begin{thebibliography}{99.}

\bibitem{A} 
Awodey, S., Sheaf representation for topoi. 
\emph{Journal of Pure and Applied Algebra}, 145, pp.~107--121, 2000.

\bibitem{AF} 
Awodey, S. and H.~Forssell, First-order logical duality. 
\emph{Annals of Pure and Applied Logic}, 164(3), pp.~319--348, 2013.

\bibitem{For}
Forssell, H., Topological representation of geometric theories. 
\emph{Mathematical Logic Quarterly}, 58, pp.~380-393, 2012.

\bibitem{AB} 
Awodey, S. and S.~Breiner, Scheme representation for first-order logic.
TACL 2013. Sixth International Conference on Topology, Algebra and Categories in Logic,
pp.~10--13, 2014.

\bibitem{B}  
Breiner, S., \emph{Scheme Representation for First-Order Logic}. 
Ph.D.~thesis, Carnegie Mellon University, 2013.  Available as {\tt arXiv:1402.2600}.

\bibitem{BM}  
Butz, C.\ and Ieke Moerdijk, Representing topoi by topological groupoids. 
\emph{Journal of Pure and Applied Algebra}, 130(3), pp.~223--235, 1998.

\bibitem{F}
Freyd, P., Aspects of topoi, 
\emph{Bulletin of the Australian Mathematical Society} 7, pp.~1--76, 1972.
 
\bibitem{H} 
Henkin, L., Completeness in the theory of types. 
\emph{Journal of Symbolic Logic} 15, pp.~81--91, 1950.

\bibitem{J}  
Johnstone, P.T., \emph{Stone Spaces}. Cambridge Studies in Advanced Mathematics 3,
Cambridge University Press, 1982.

\bibitem{JM}  Joyal, A.\ and I.~Moerdijk, Toposes as homotopy groupoids. 
\emph{Advances in Mathematics}, 80(1), pp.~22--38, 1990.

\bibitem{JT}  Joyal, A.\  and M.\ Tierney, An extension of the Galois theory of Grothendieck.
\emph{Memoirs of the AMS}, 308, 1984.

\bibitem{L1} 
Lambek, J., On the sheaf of possible worlds. 
In Adamek, J.\ and  Mac Lane, S.\ (Eds.), \emph{Categorical Topology}, World Scientific, Singapore, 1989.

\bibitem{L2} Lambek, J., What is the world of mathematics? 
\emph{Annals of Pure and Applied Logic}, 126, pp.~149--158, 2004.

\bibitem{LM}  
Lambek, J. and I.~Moerdijk, Two sheaf representations of elementary toposes. 
In A.S.\ Troelstra and D.\ van Dalen (Eds.), \emph{Brouwer Centenary Symposium}, North-Holland, Amsterdam, 1982.
 
\bibitem{M} 
 Makkai, M., Stone duality for first order logic. 
 \emph{Advances in Mathematics}, 65(2), pp.~97--170,  1987.
 
 \bibitem{MM} 
Mac~Lane, S.\ and I.~Moerdijk, \emph{Sheaves in Geometry and Logic}. Universitext. Springer, 2nd edition, 1992.

\end{thebibliography}
\end{document}